\documentclass[12pt]{article}

\usepackage{amssymb, amsthm, amsgen, amscd, amssymb, amsmath, amsxtra, subcaption, caption}

\usepackage{graphics, epsfig, color}

\usepackage[parfill]{parskip}

\usepackage[matrix, arrow, curve]{xy}

\usepackage{amsmath,amscd}

\usepackage[pdftex, plainpages=false, pdfpagelabels]{hyperref}

\usepackage[letterpaper, margin=1in]{geometry}

\hypersetup{
    linktocpage=true,
    colorlinks=true,
    bookmarks=true,
    citecolor=blue,
    urlcolor=blue,
    linkcolor=blue,
    citebordercolor={1 0 0},
    urlbordercolor={1 0 0},
    linkbordercolor={.7 .8 .8},
    breaklinks=true,
    pdfpagelabels=true,
    }

\theoremstyle{remark}

\def\bee{\begin{equation*}}

\def\eee{\end{equation*}}

\def\be{\begin{equation}}

\def\ee{\end{equation}}

\def\R{{\mathbb R}}

\definecolor{Maroon}{rgb}{0.8, 0.0, 0.0}

\title{\vskip-1in Folklore, the Borromean rings, the icosahedron, and three dimensions}
\author{Dave Auckly}

\date{}

\begin{document}
\bibliographystyle{plain}
\maketitle
\vskip-.2in
 \centerline{Mathematics Department, Kansas State University,
Manhattan, KS 66506 }
\centerline{{\tt{dav@math.ksu.edu} }}

\section*{Introduction}\label{intro}

There are many bits of mathematical folklore --  results that have been known for a long time -- that are not covered in a typical undergraduate curriculum. We will give a proof that the Borromean rings are linked, but our proof will not be the shortest possible. Indeed, we give a number of unnecessary, but loosely related, stories to expose more folklore. In particular we will see that the Borromean rings are related to the icosahedron and something called the Poincar\'e  
homology sphere. We begin with a cautionary tale.

Around 1900 Poincar\'e made many important contributions to topology. He was considering spaces that looked a bit like $3$-space, but were not necessarily just $\R^3$. The same phenomena may be seen with shapes that look like $2$-space. For instance a small neighborhood on the surface of the earth looks like $\R^2$ even though the surface of the Earth is not just a flat plane. The surface of the Earth is better modeled by the $2$-sphere $S^2:=\{v\in\R^3\,|\, |v|^2=1\}$. The same could be said about the surface of a donut, neighborhoods of it look like $\R^2$ but it isn't. Shapes that look like $\R^2$ are called $2$-manifolds and shapes that look like $\R^3$ are called $3$-manifolds.  

Poincar\'e noticed that a compact, roughly speaking ``does not extend to infinity," $2$-manifold that has no holes in it must be the two-sphere, $S^2$. He conjectured that the same thing was true for $3$-manifolds where the measure of holes in a manifold $M$ is given by something that we now call the first homology, $H_1(M)$.  In 1904 he found a counter-example  to this conjecture. To see this he introduced the notion of the fundamental group of a space, $\pi_1(M)$, and constructed a counter-example called the Poincar\'e  homology sphere that may be described using an icosahedron. He then revised his conjecture to state that the only compact $3$-manifold with trivial fundamental group is $S^3$ \cite{49still}. More than one hundred years later Perelman gave the first correct proof of this conjecture \cite{P}.

In fact, before the correct proof a number of incorrect proofs were given of the Poincar\'e conjecture.
In many incorrect proofs the only hypothesis that is used is that the homology is trivial. Of course any such proof is doomed to fail as the Poincar\'e  homology sphere shows. John Stallings wrote a nice paper about this called, {\it How Not to Prove the Poincar\'e  Conjecture} \cite{stal}. Had these authors payed a bit more attention to the folklore about the difference between homology and fundamental group, there may have been fewer false proofs. The  first homology $H_1(M)$ is the abelianization of the fundamental group, $\pi_1(M)$. The difference between trivial homology and trivial fundamental group, while understood, can be subtle. We will explain the fundamental group in greater detail in the first section below.

We note that the electronic version of this paper has color graphics that may be easier to follow.

\section*{The fundamental group} 
\begin{figure}[!ht]  
\centerline{
\includegraphics[width=95mm]{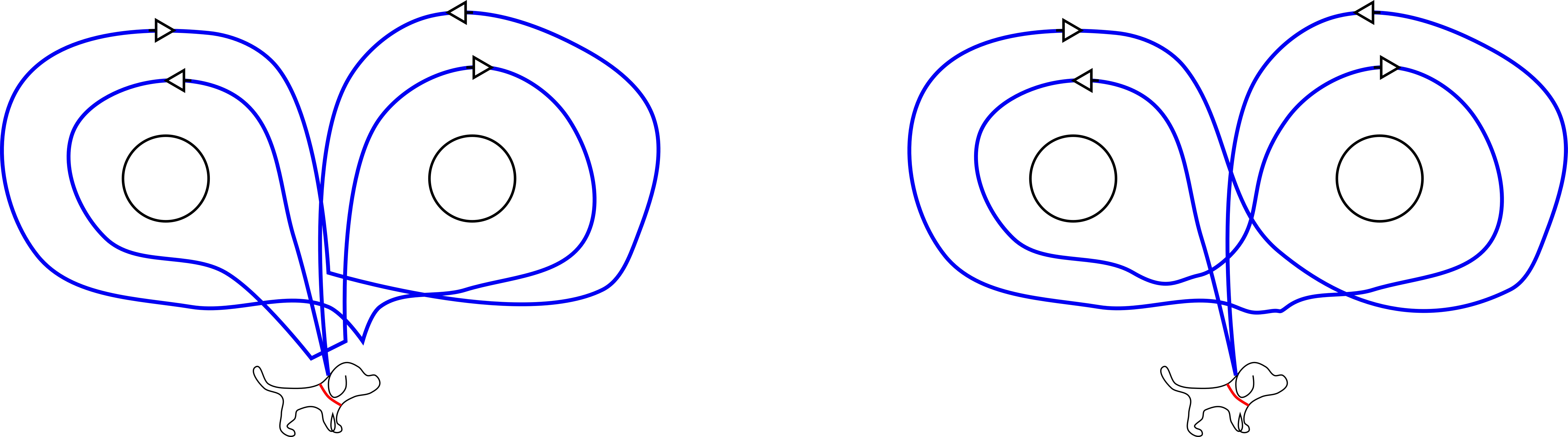}}
\caption{Tangled leash}\label{dog-walk}
\end{figure}

Imagine you are walking a slightly skittish dog on a leash in a forest. Most of the time the dog will stay by your side. Now the dog thinks it detects a squirrel and sniffs its way around a tree twice. You are stuck until you can convince the dog to unwind from the tree. There is a group here. The dog could sniff $2$ times around the tree, or $13$ times, or could go the other way to wrap $-7$ times around the tree. The amount the dog winds around the tree can be measured by an integer.  The leash tangled around two trees is displayed in Figure~\ref{dog-walk} on the next page.

A more precise description of this may be given using the fundamental group. Let $X$ be a topological space, i.e., a space for which the notion of continuous is defined, and let $x_0$ be a point in $X$. Elements of the fundamental group are equivalence classes of continuous paths,
$\gamma:[0,1]\to X$ with $\gamma(0)=\gamma(1)=x_0$. Two paths, $\gamma_0$ and $\gamma_1$ are equivalent 
if there is a deformation (\emph{homotopy}) $H:[0,1]\times [0,1] \to X$ so that $H(0,s)=\gamma_0(s)$, $H(1,s)=\gamma_1(s)$, $H(t,0)=x_0$, and $H(t,1)=x_0$. In the example of the dog walk, the space $X$ is all of the ground that is not covered by the tree. The base point $x_0$ is the point that you are standing on. The leash may be viewed as a continuous path $\gamma$. In the equivalence relation, one can view the parameter $t$ in the function $H$ as time. Thus   
 $H(0,s)=\gamma_0(s)$ requires that the leash is along the path $\gamma_0$ at time zero. Fixing $t$, the function $H(t,\cdot)$ represents the position of the leash at time $t$. The condition $H(t,0)=x_0$ corresponds to the handle of the leash being with you at all time, and the condition $H(t,1)=x_0$ corresponds to your skittish dog sitting on your foot for all time. This is captured in the definition below.
\[
\pi_1(X,x_0) := \{\gamma : [0,1]\to X\,|\, \gamma(0)=\gamma(1)=x_0 \}/\text{homotopy}\,.
\]

The group operation is called \emph{concatenation}. Let $\gamma$ represent the position of the leash after one short exploration of the dog, and let $\delta$ represent the position of the leash after another short exploration. The product, $\gamma *\delta$ is what would happen if the dog first did $\gamma$, and then did $\delta$. The formula is:
\[
(\gamma *\delta)(t)\,=\,\begin{cases}
\gamma(2t) & \text{if} \ t\in [0,1/2],\\
\delta(2t-1) & \text{if} \ t\in [1/2,1].
\end{cases}
\]

Going further on this walk you now approach a pair of close trees, and there is a squirrel in them. In complete excitement your dog races counterclockwise around the left tree $L$ back to you, then clockwise around the tree on the right $R^{-1}$ and back to you, clockwise around the tree on the left $L^{-1}$, and counterclockwise around the tree on the right $R$. A figure showing the leash and how it may be pulled a bit tighter to an equivalent (homotopic) path is shown below. This path will be denoted by $LR^{-1}L^{-1}R=L*R^{-1}*L^{-1}*R$. 

There is something interesting about this particular example. In total the dog has wrapped around the left tree zero times, so the leash would pull free if the right tree was not there. In addition the dog walked around the right tree a total of zero times. Yet as long as the squirrel continues to scare your dog and your dog does not leave your foot, the leash will be linked with the trees. This is exactly the difference between the first homology and the fundamental group. The first homology is the abelianization of the fundamental group (in other words $XY=YX$ in homology). In homology $LR^{-1}L^{-1}R=LL^{-1}R^{-1}R=1$. We typically write the operation in the fundamental group multiplicatively ($XY$) and operation in the homology additively ($X+Y$) to remind us that the homology is abelian. 

The fundamental group of the complement of two disjoint disks in the plane (the forest floor away from the two trees) is what is known as a \emph{free group} on two generators. The elements of this group are a trivial element $1$ and all finite ``words" that may be bade from the letters $L$, $R$, $L^{-1}$, $R^{-1}$ with all cancelations of $X$ next to $X^{-1}$. Multiplication is just the process used to create a compound word. Thus the product of $LRLR^{-1}$ with $RRRLL$ is $LRLRRLL$.

The false proofs of  the Poincar\'e conjecture asserted that when every loop in a compact $3$-manifold wrapped zero times around any hole an the abelian sense as measured in homology, the manifold had to be trivial. This turned out to be false. Now that we have a new trick for our dog, let's find a better place to tangle our leash.

\section*{The Borromean rings}
\begin{figure}[!ht]  
\centering
\vspace{-12pt}
    \begin{subfigure}[b]{0.45\textwidth}
\includegraphics[width=\textwidth]{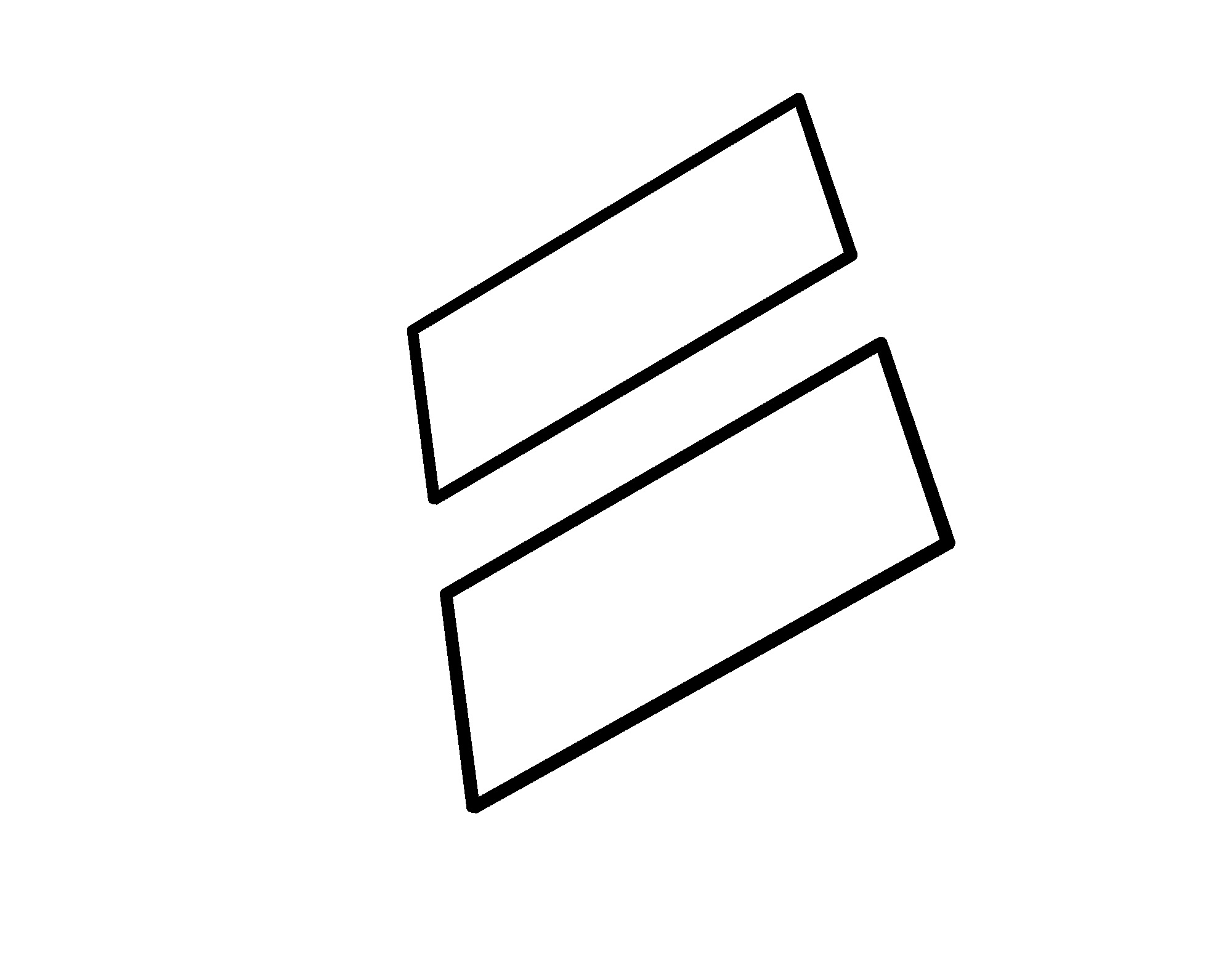}
\caption{Loose}\label{2loose} 
\end{subfigure}
\begin{subfigure}[b]{0.4\textwidth}
\includegraphics[width=\textwidth]{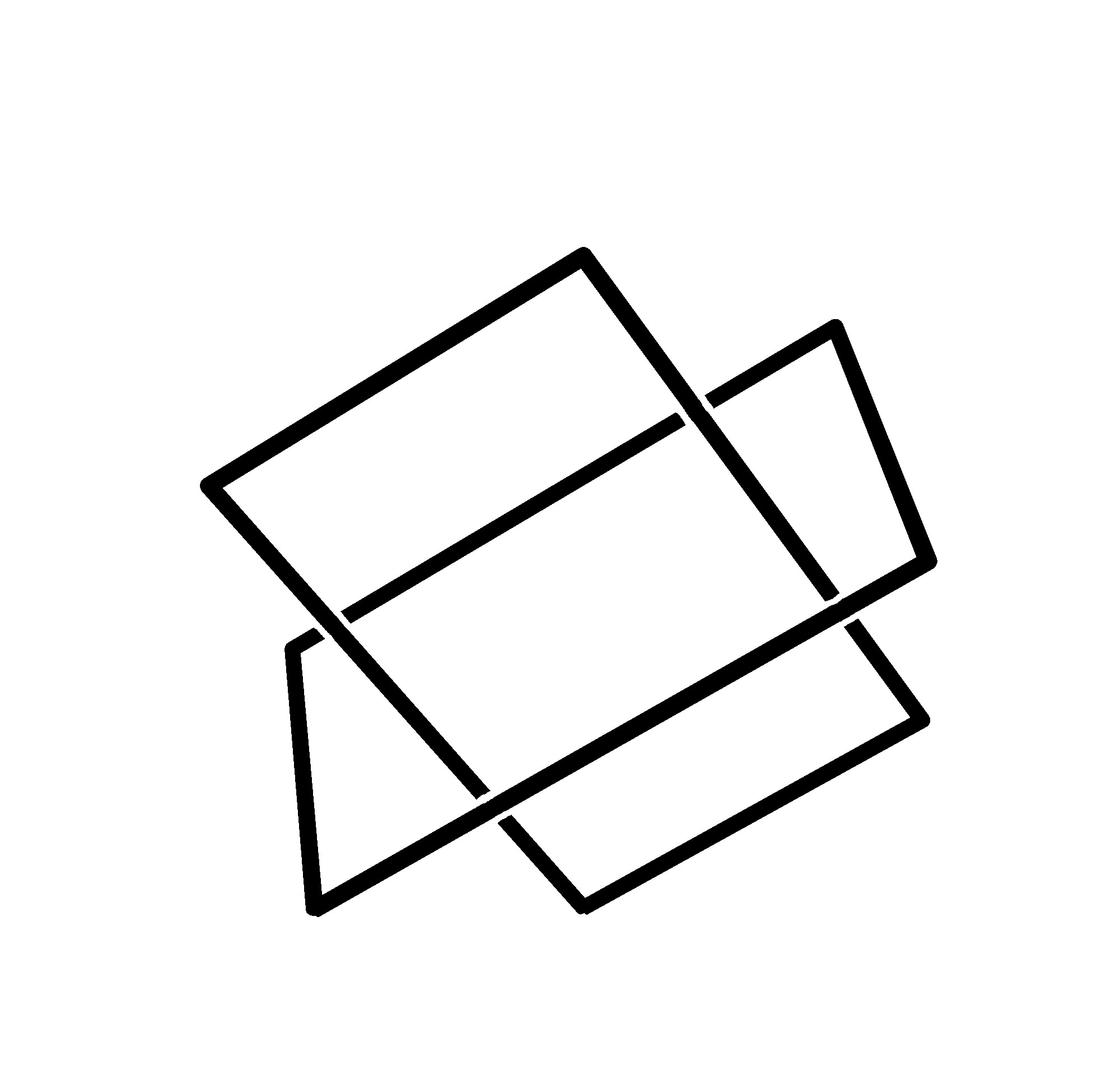}
\caption{In position}\label{2pos}
\end{subfigure}
\caption{Two unlinked rectangles}\label{2rect} 
\end{figure}

Start by making a roughly $34$ inch by $55$ inch rectangle. Really a $34$ by $34\cdot(1+\sqrt{5})/2$ rectangle. The number $(1+\sqrt{5})/2$ is the \emph{golden ratio}. It is the length of the diagonal of a sidelength $1$ pentagon. These numbers are close to adjacent numbers in the Fibonacci sequence ($34$ and $55$) and this is a nice size for public display. We could prepare and make two rectangles as in Figure \ref{2loose} and later move them into position as in  Figure \ref{2pos}.

\begin{figure}[!ht]  
\centering
\vspace{-12pt}
    \begin{subfigure}[b]{0.45\textwidth}
\includegraphics[width=\textwidth]{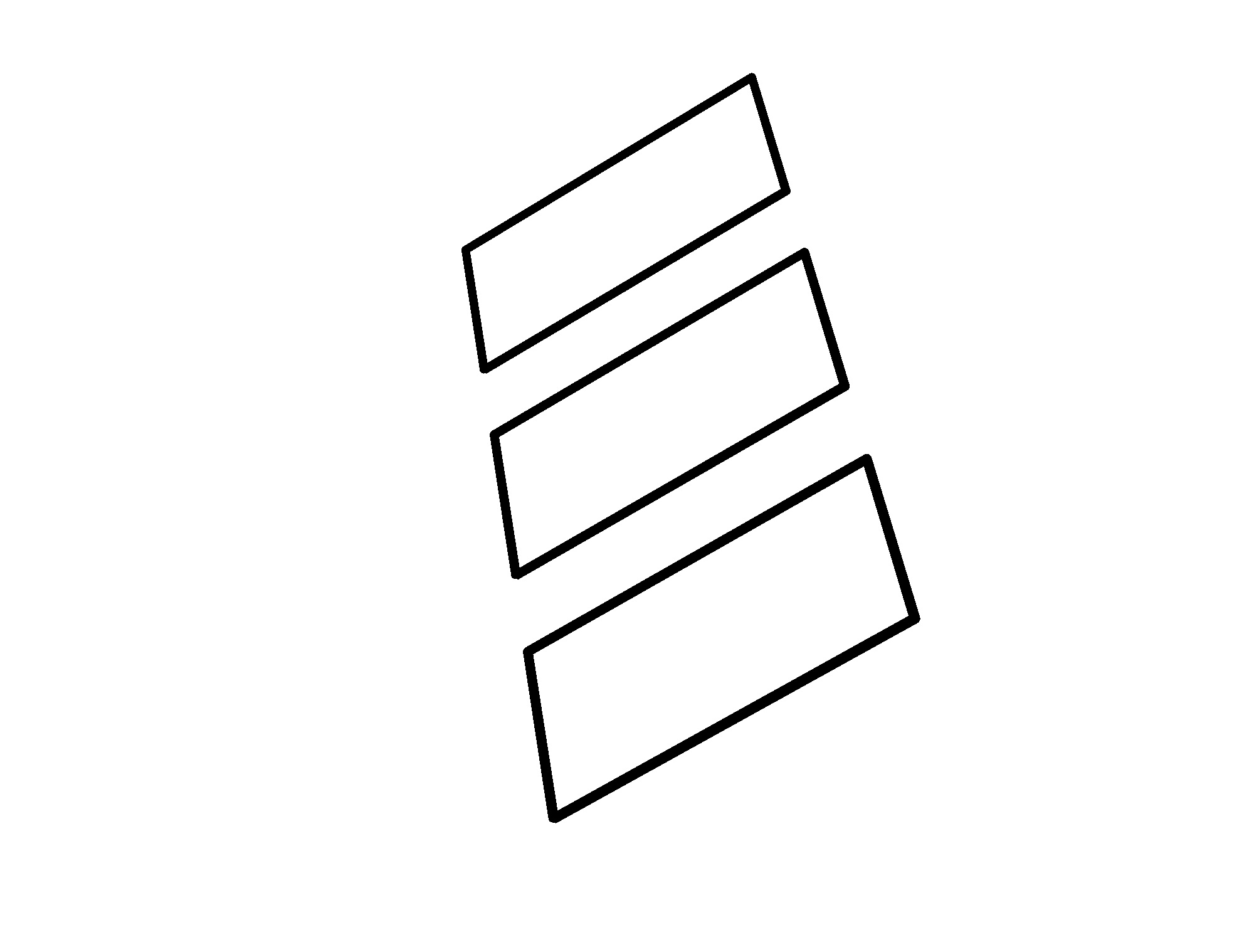}
\caption{Unlinked}\label{3a}
\end{subfigure}
\begin{subfigure}[b]{0.4\textwidth}
\includegraphics[width=\textwidth]{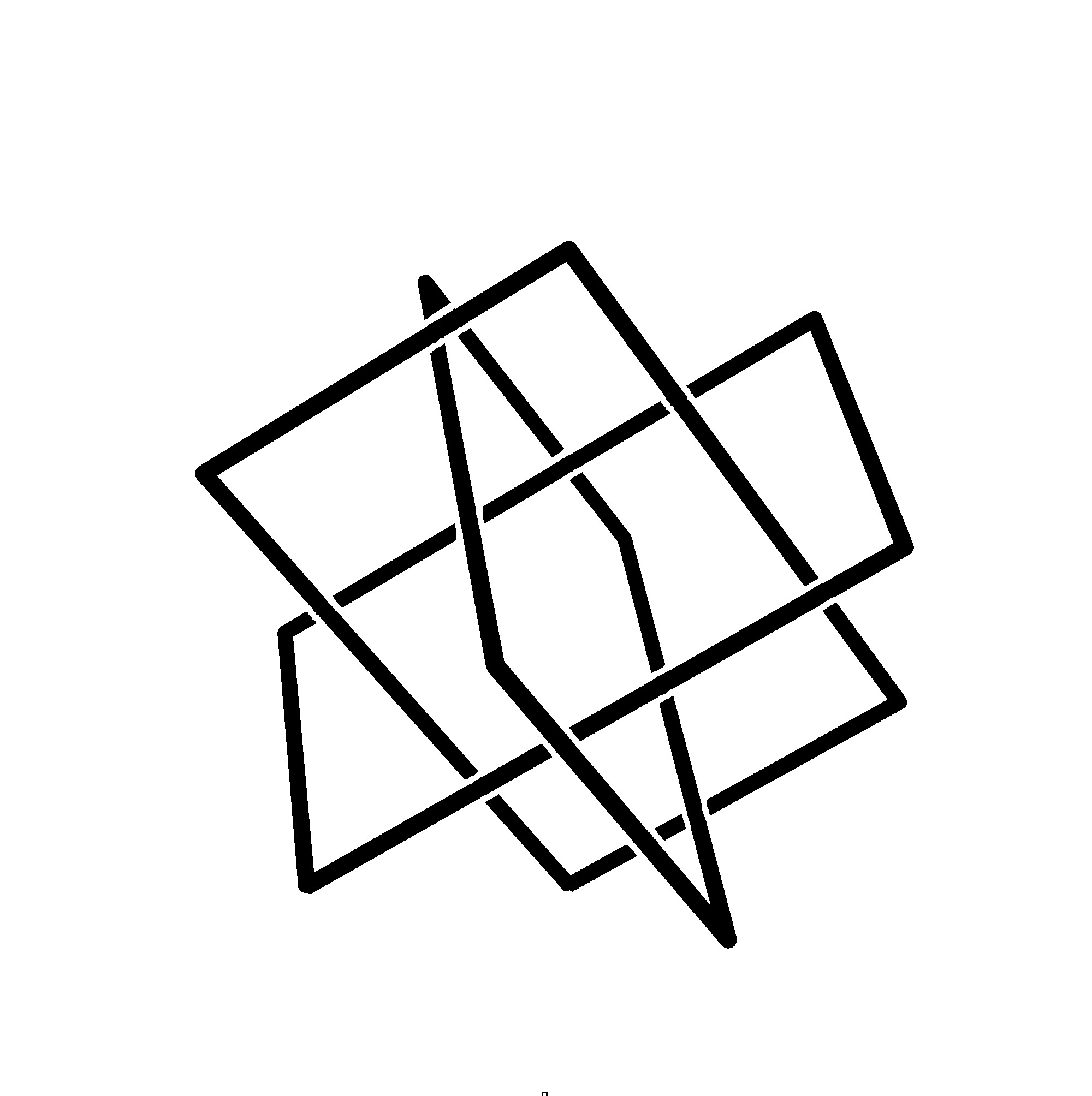}
\caption{Linked?}\label{3b}
\end{subfigure}
\caption{Assemble three rectangles}\label{3rect}
\end{figure}

\begin{figure}[!ht]  
\centering
\vspace{-1.2pt}
    \begin{subfigure}[b]{0.45\textwidth}
\includegraphics[width=\textwidth]{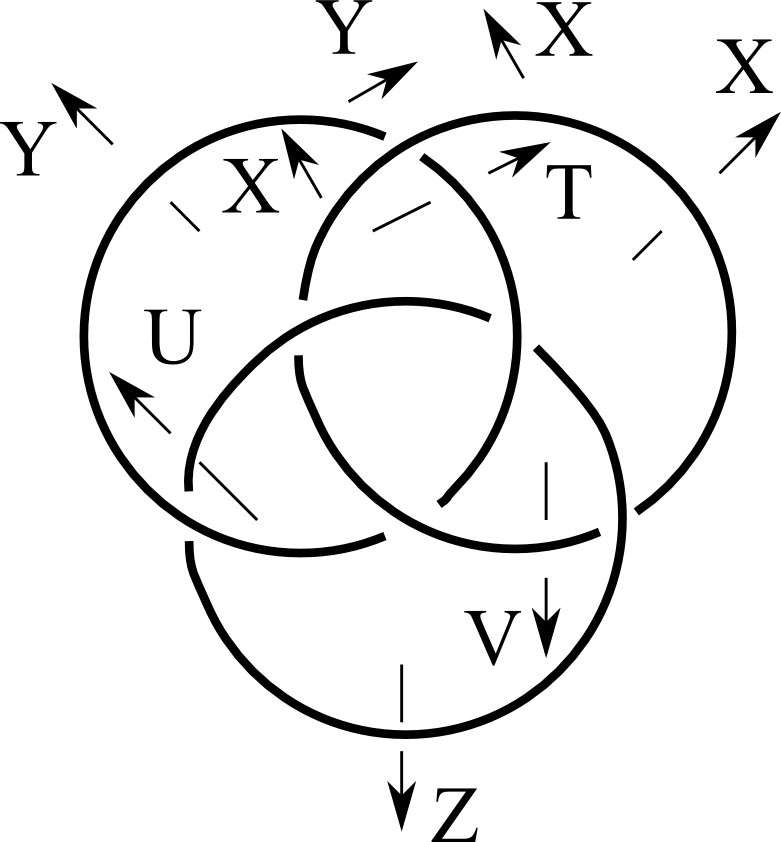}
\caption{Generators}\label{5a}
\end{subfigure}
\begin{subfigure}[b]{0.45\textwidth}
\includegraphics[width=\textwidth]{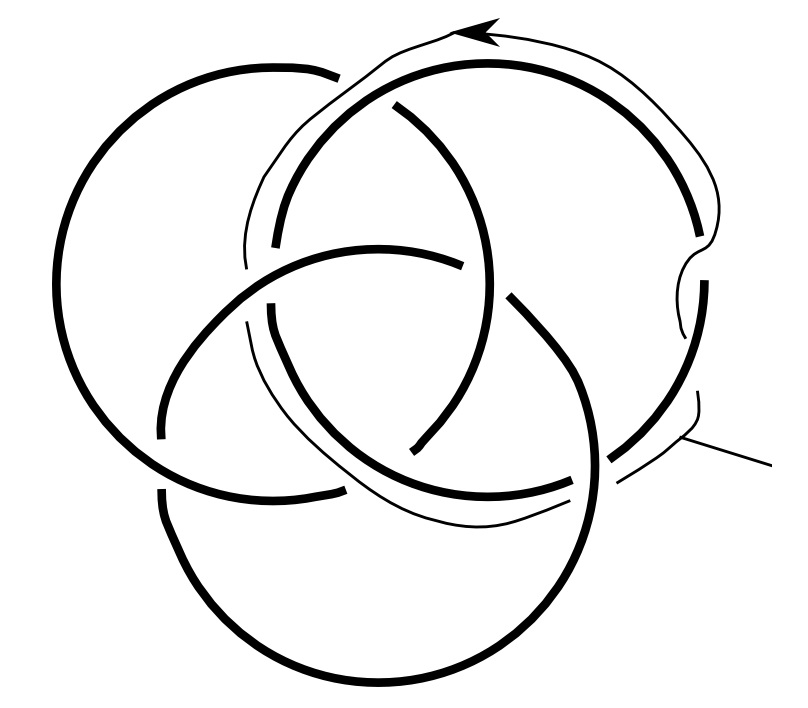}
\caption{Glue curve}\label{5b}
\end{subfigure}
\caption{Borromean Rings.}\label{bor}
\end{figure}

A more precise description of the right hand side starts with three of these of these rectangles. We put one with center at the origin and line of symmetry parallel to the long side (call this the rectangle axis) matching the $x$-axis of the $x-y$ plane. Similarly, we could put one with the axis matching the $y$-axis of the $y-z$ plane and one with the axis matching the $z$-axis of the $z-x$ plane. The result is the very specific model of the  Borromean rings in Figure \ref{3b}. The Borromeo family used a version of this link on its coat of arms. The version from the coat of arms is the traditional image of the Borromean rings, and it is displayed in Figure \ref{bor} later.

This brings us to a question:
\begin{quote}
{\bf Question:} If we set our three rectangles out on the floor as in  Figure \ref{3a}, would we be able to move them into the 
desired configuration in Figure \ref{3b} without taking one apart?
\end{quote}

The answer is that it is not possible to assemble three rectangles into the Borromean configuration without breaking one of the rectangles. We will ultimately give a proof of this fact using the fundamental group of the complement of the rectangles.

\begin{figure}[!ht]  
\centering
\vspace{-1.2pt}
    \begin{subfigure}[b]{0.3\textwidth}
\includegraphics[width=.7\textwidth]{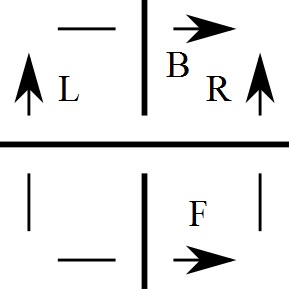}
\caption{Top view}\label{4a}
\end{subfigure}
\begin{subfigure}[b]{0.3\textwidth}
\includegraphics[width=\textwidth]{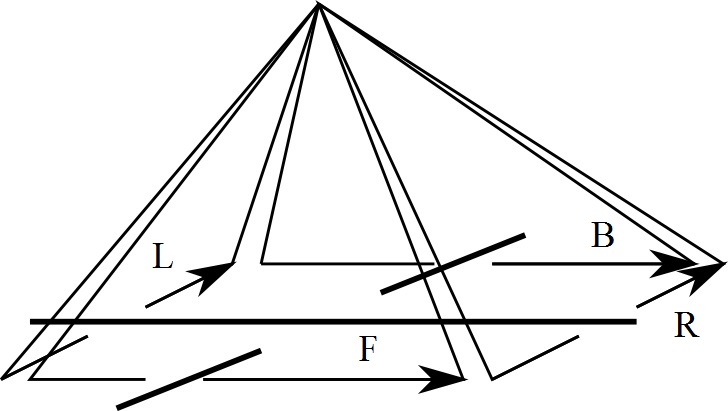}
\caption{Super leash}\label{4b}
\end{subfigure}
\begin{subfigure}[b]{0.3\textwidth}
\includegraphics[width=\textwidth]{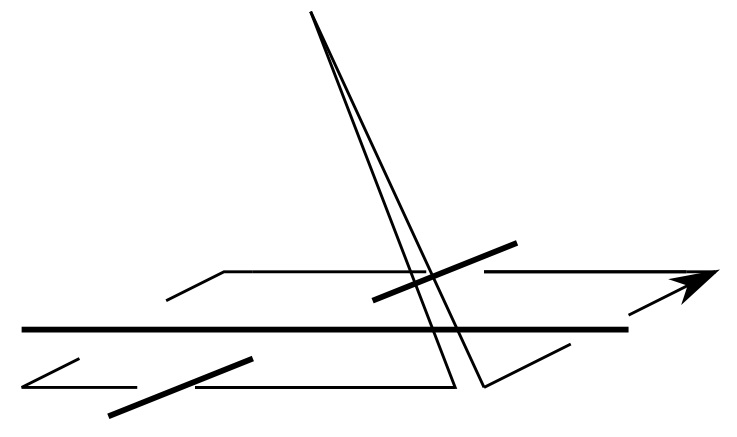}
\caption{Deformed}\label{4c}
\end{subfigure}
\caption{Loops Near a Crossing}\label{wgen}
\end{figure}

A \emph{link} is a disjoint collection of topological circles in $3$-space.
Given a picture of a  link, one may label some loops by putting arrows crossing under various strands of the link as in the Figure \ref{4a}. We just discovered that your dog is a super dog. Your dog now flies up to a point above a crossing.  Each arrow in Figure \ref{wgen} represents the loop of leash generated when your dog starts in the air (at the base point) flies to the tail of the arrow, follows the arrow and returns to the floating base point. The arrows labeled $L$ and $R$ are below the horizontal strand but above the height of the vertical strand in Figure \ref{4a}.Thus the loops represented by arrows labeled $L$ and $R$ are equivalent because one may be deformed into the other. However, the loops labeled $F$ and $B$ are not equivalent, because the horizontal component of the link gets in the way of the deformation. See Figure \ref{4b}. The inverse of a loop, is the same loop, with reversed orientation.

Superdog will now show us a trick while trying to tie its leash around the strands. Superdog flies
\begin{enumerate}
\item down then along the arrow on the right,
\item back up, back down, 
\item follows the arrow in the back in the opposite direction, 
\item back up, then down, 
\item follows the arrow on the left in the opposite direction, 
\item back up, then down, 
\item follows the front arrow in the indicated direction, and 
\item finally fies back up to the base point. 
\end{enumerate} 
The result is 
the product $RB^{-1}L^{-1}F$ displayed in Figure \ref{4b}. However, this is the trivial element of the fundamental group. To see this, note that the portions of the leash generated by the "back up, then down" portions 2, 4, and 6 of the flight may be deformed down to give the representative in Figure \ref{4c}. From here the square portion of the leash loop on the bottom can shrink and the entire loop can be pulled back up to the base point.
It is a fact that the fundamental group of the complement of any link may be described as the group with one generator for each strand and one relation similar to the one we just described for each crossing. This is called the  \emph{Wirtinger presentation}. See the book by Dale Rolfsen for a proof of this, \cite{R}.

How could this possibly help us? We are trying to show that we can not move the unlinked rectangles from Figure \ref{3a} into the configuration in Figure \ref{3b}. Imagine that we could, and fill up the space around the rectangles with honey. When one pushes the set of rectangles from Figure \ref{3a} over to the position in Figure \ref{3b}, all of the honey will move as well. The correspondence between the starting location of a molecule of honey and the ending location of that molecule of honey gives a homeomorphism  between the complement of the configuration in Figure \ref{3a} and the complement of the configuration in Figure \ref{3b}. Thus the fundamental groups of the complements of these two configurations would be the same. The fancy name for this idea is the isotopy extension theorem.

The fundamental group of the complement of the trivial $3$-component link in Figure \ref{3a} is the group with three generators and no relations (one generator for each component, and no relations because there are no crossings):
\[
\pi_1(\text{Complement of unlink}) = \langle P, Q, R \rangle\,.
\]

The three rectangles in the configuration in Figure \ref{3b} may be continuously deformed, without crossing, to the configuration in  Figure \ref{5a}. Thus the fundamental groups of the complements of these configurations are the same.




To compute the fundamental group of the complement of the Borromean rings we label arrows under the strands in Figure \ref{5a} with $T,U,V,X,Y,Z$. The strand on the upper right of this figure has an $X$ arrow running under it. The loop represented by this arrow starts at a base point above the figure (say your nose), travels down to the tail of the arrow, follows the arrow and returns to the basepoint. This loop may be slid (homotoped) to other positions along this strand as long as one does not try to slide past an undercrossing. We label two loops equivalent to $X$ near the top crossing to make it easier to see the relation arising from this crossing. Following the procedure in Figure \ref{wgen} starting at the head of the $T$ and prodeeding clockwise we read the relation 
$T^{-1}XYX^{-1}=1$.
Thus $T=XYX^{-1}$. Similarly, the outside crossing on the left shows that $U=YZY^{-1}$ and the outside crossing on the right shows that $V=ZXZ^{-1}$. Using the three inner crossings leads to the presentation:
\[
\begin{aligned}
\pi_1(\text{Complement of the Borromean rings}) = &\\
 \langle X, Y, Z \,  |\,  YZY^{-1}XYZ^{-1}Y^{-1}ZX^{-1}Z^{-1} &=1, \\
 XYX^{-1}ZXY^{-1}X^{-1}YZ^{-1}Y^{-1} &=1, \\
 ZXZ^{-1}YZX^{-1}Z^{-1}XY^{-1}X^{-1} &=1 \rangle\,.
\end{aligned}
\]
Indeed, the relation coming from the inner crossing closest to arrow $Y$  going counter-clockwise starting with the $U$ on the top reads: \hfill\newline \[1=UXU^{-1}V^{-1}=(YZY^{-1})X(YZ^{-1}Y^{-1})(ZX^{-1}Z^{-1})\,.\]

To prove that the two links are different, we could just show that these two groups are different.
Instead we will use a slightly weird argument. We choose a weird argument because it will allow us to discuss the relationship between the Borromean rings, the Poincar\'e homology sphere and the icosahedron. Sometimes the long way around is more scenic.

\section*{From the Borromean Rings to the Poincar\'e Sphere}

This is the most technical section of this paper. A quick summary is that if we add the relation $P=1$ and the corresponding relations for the other components ($Q=R=1$) to the fundamental group of the complement of the trivial $3$-component link, we would get a trivial group. If the trivial $3$-component link was equivalent (isotopic) to the Borromean rings the relation corresponding to $P=1$ would be $X^{-1}U^{-1}Z=1$. Adding this to the fundamental group of the complement of the Borromean rings together with the analogous relations for the other two components results in a non-trivial group so the two links are not equivalent.  Right now you should probably have no idea where the relation  $X^{-1}U^{-1}Z=1$ comes from. We will turn this around and start with the relation $X^{-1}U^{-1}Z=1$ and see that it leads to $P=1$. We will also construct a space with the corresponding fundamental group -- the Poincar\'e homology sphere.

We can perform something called surgery on the trivial $3$-component link and on the Borromean rings. If these two links were the same, the result of the corresponding surgeries would be the same. We begin with a warm-up. Consider the two ``links" in Figure \ref{pts} below. Each is two points put into a space consisting of two $2$-spheres.

\begin{figure}[!ht]  
\centerline{
\includegraphics[width=109mm]{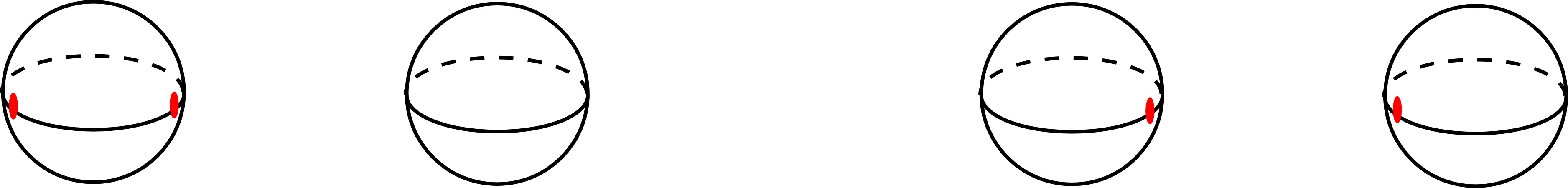}}
\caption{The points}\label{pts}
\end{figure}

It there was a deformation taking the trivial $2$-point link on the left to the $2$-point link on the right, there would be a homeomorphism taking the complement of the two points on the left to the complement of the two points on the right. (This is the isotopy extension theorem again. Fill the complement with honey and see where the honey molecules would go after pushing the one link to the other.)

We will now glue the same thing to each complement, namely a cylinder. Notice that $S^1\times (0,1)$ is homeomorphic to the open unit ball with the origin deleted. Indeed, just consider the coordinate in $S^1$ as an angle and the coordinate in $(0,1)$ as a radius and   $S^1\times (0,1)$ will be the punctured open ball expressed in polar coordinates. Similarly $S^1\times (-1,0)$ is also homeomorphic to a deleted open ball. Thus we can glue the open cylinder, $S^1\times (-1,1)$ to the complement of either link. The homeomorphism that we get by assuming that the links are equivalent extends and would imply that the result of gluing in the open cylinder to the left side would be homeomorphic th the result of gluing the open cylinder to the right side. However, we display the result of this gluing in Figure \ref{S0sur} below, and can tell that the two sides are not homeomorphic because one is disconnected, and the other is connected.
This process is called surgery because we are cutting something out and sewing something else back in, a bit like an organ transplant.

\begin{figure}[!ht]  
\centerline{
\includegraphics[width=109mm]{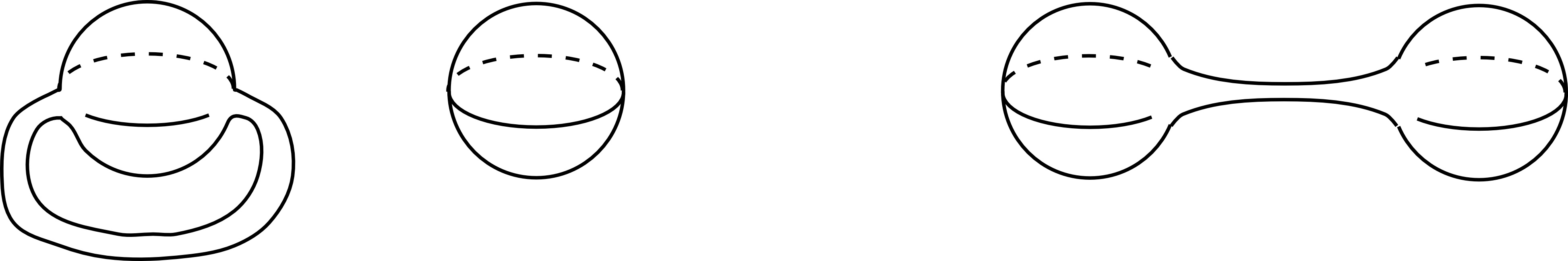}}
\caption{Surgery on the points}\label{S0sur}
\end{figure}

We are now going to do the analogous thing to the Borromean rings.We cut out the rings and glue in solid tori (donuts). This is displayed in Figure \ref{poin} below.

\begin{figure}[!ht]  
\centering
\vspace{-1.2pt}
    \begin{subfigure}[b]{0.45\textwidth}
\includegraphics[width=\textwidth]{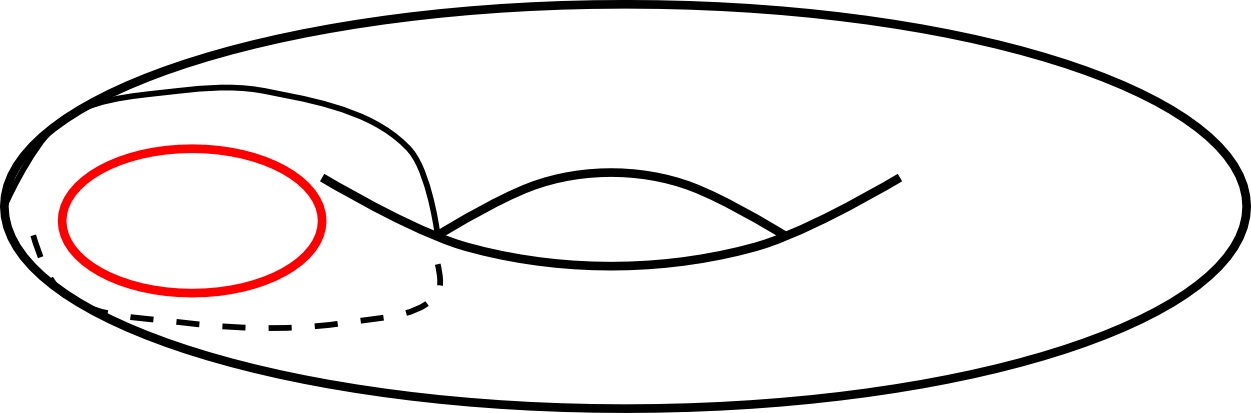}
\caption{Donut}\label{Pa}
\end{subfigure}
\begin{subfigure}[b]{0.45\textwidth}
\includegraphics[width=\textwidth]{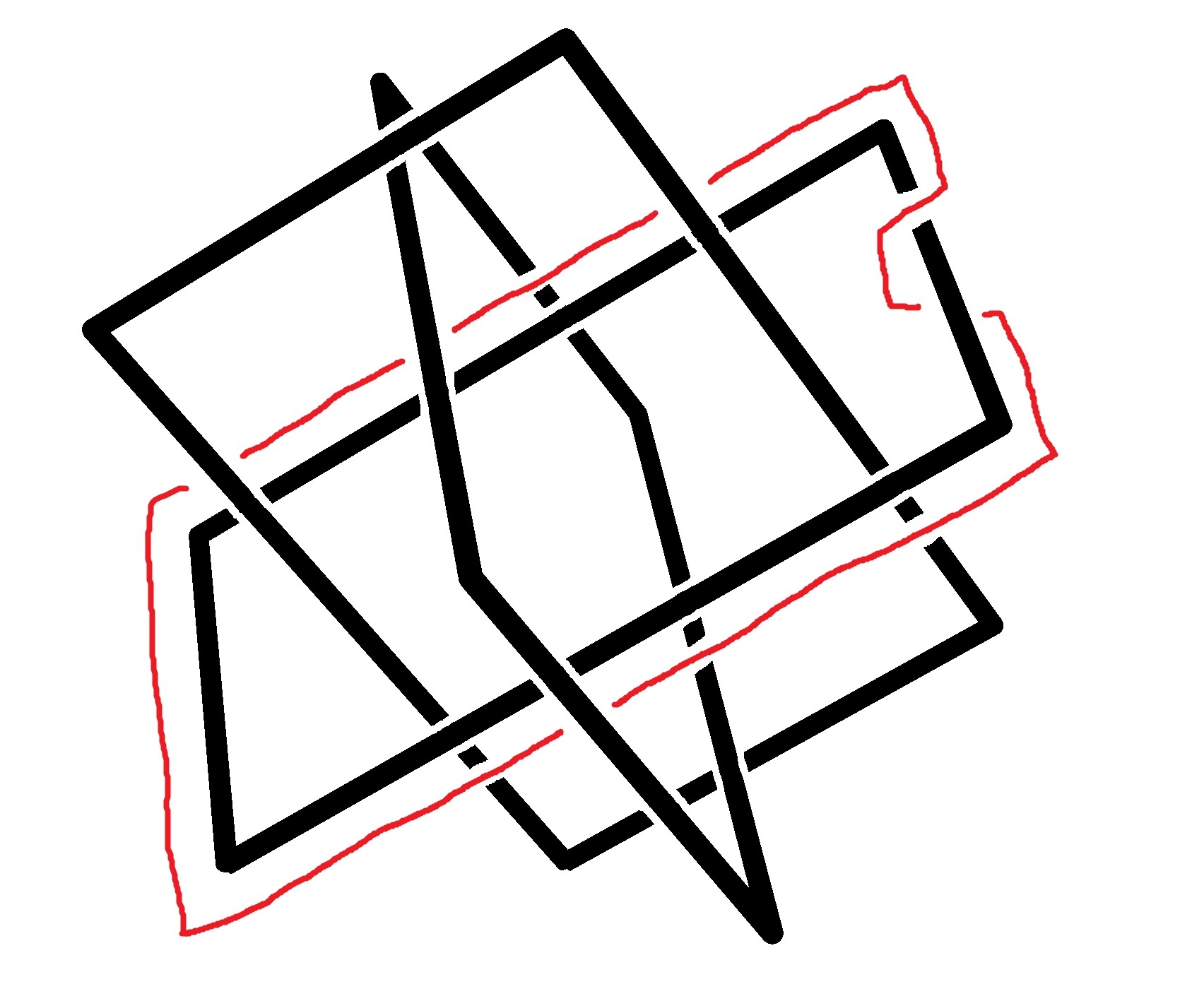}
\caption{Rectangular glue curve}\label{Pb}
\end{subfigure}
\caption{The Poincar\'e homology sphere.}\label{poin}
\end{figure}

To specify how an open torus is glued to the complement of the Borromean rings, we will keep track of one curve. The thin black curve on the left side of the boundary of the solid torus on the left of Figure \ref{poin} bounds a disk. However it is not part of the open solid torus.
 To be clear
the open solid torus is given by 
\[
\{(v,w)\in \R^2\times \R^2\,|\, |V|<1, \ |w|=1\}\,.
\]
The red curve is parallel to the thin black curve. When we glue the open solid torus to the complement of the Borromean rings we will make sure that the thin (red in electronic) curves in Figure \ref{poin} match. The thin red curve in the complement of the Borromean rings is called a surgery curve. The reason we keep track of this curve is because the corresponding loop will be trivial in the new manifold. We do the same gluing with each of the other components. 

\begin{quote}
The result of attaching three open solid tori to the complement of the Borromean rings is equivalent to removing one point from the Poincar\'e homology sphere. We can take this as the definition of the Poincar\'e homology sphere.
\end{quote}

\begin{figure}[!ht]  
\centerline{
\includegraphics[width=39mm,angle=90,origin=c]{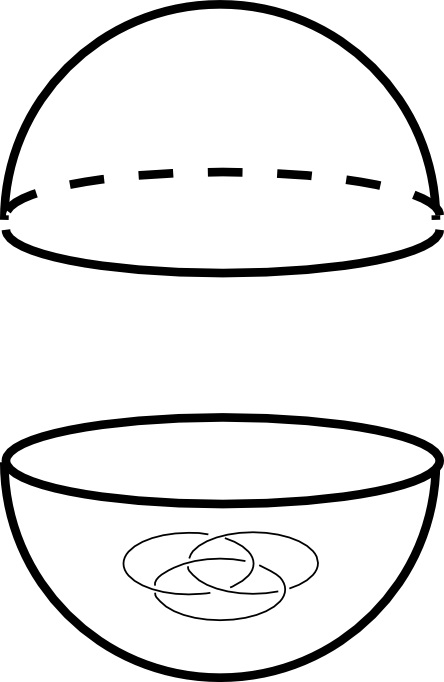}}
\caption{Two balls}\label{cap}
\end{figure}

We could fill in the missing point by gluing in a hemisphere from the $3$-dimensional sphere, i.e. the set of points one unit from the origin in $\R^4$. A schematic of this is displayed in Figure \ref{cap}. Here the left hemisphere with a small copy of the
Borromean rings represents the  result of removing the Borromean rings from a solid ball and gluing in donuts. The right hemisphere represents a second solid ball. Gluing these together yields the  Poincar\'e homology sphere. We denote it by $\Sigma$. To compute the fundamental group of Poincar\'e homology sphere we just need to add relations indicating that the loop corresponding to the thin (red) circle and the analogous loops are trivial. The thin (red) loop  in Figure \ref{poin} corresponds to the loop in the  Figure \ref{5b}. Reading the relation we  get $X^{-1}U^{-1}Z=X^{-1}YZ^{-1}Y^{-1}Z=1$.
 Adding the analogous relations for the other two components gives:
\[
\begin{aligned}
\pi_1(\Sigma) = \langle X, Y, Z \, |\, & YZY^{-1}XYZ^{-1}Y^{-1}ZX^{-1}Z^{-1}=1, \\ & XYX^{-1}ZXY^{-1}X^{-1}YZ^{-1}Y^{-1}=1, \\
& ZXZ^{-1}YZX^{-1}Z^{-1}XY^{-1}X^{-1}=1, \\
& X^{-1}YZ^{-1}Y^{-1}Z=1, \\
& Y^{-1}ZX^{-1}Z^{-1}X=1, \\
&Z^{-1}XY^{-1}X^{-1}Y=1  \rangle\,.
\end{aligned}
\]
This is the fundamental group of the Poincar\'e homology sphere. We can simplify the expression for this group.
The last relation in $\pi_1(\Sigma)$ implies $Z=XY^{-1}X^{-1}Y$. Substituting
this into the other relations allows us to write the group without using the generator $Z$. Using this, the relation $Y^{-1}ZX^{-1}Z^{-1}X=1$ is seen to be equivalent to $XY^{-1}X^{-1}YX^{-1}Y^{-1}X=1$, and the relation $X^{-1}YZ^{-1}Y^{-1}Z=1$ is seen to be equivalent to  $YX^{-1}Y^{-1}XY^{-1}X^{-1}Y=1$. The other relations all follow from these two and the expression for $Z$. Thus,
\[
\pi_1(\Sigma) = \langle X, Y \, |\, XY^{-1}X^{-1}YX^{-1}Y^{-1}X=1, YX^{-1}Y^{-1}XY^{-1}X^{-1}Y=1  \rangle\,.
\]

To write this in an interesting way set $X = A^{-1}BA^{-2}B$ and $Y = A^{-1}B$. Notice that this implies that $A = YX^{-1}Y$ and $B = YX^{-1}Y^2$ so we can either use generators $X$ and $Y$ to describe the group or generators $A$ and $B$. Using this substitution the relation $YX^{-1}Y^{-1}XY^{-1}X^{-1}Y = 1$ becomes $AB^{-2}A^2 = 1$ or just $B^2 = A^3$. Now $XY^{-1}X^{-1}YX^{-1}Y^{-1}X = 1$ becomes $A^{-1}BA^{-1}B^{-1}A^2B^{-1}A^2B^{-1}A^{-1}B = 1$. Using $A^2 = B^2A^{-1}$ this becomes $(A^{-1}B)^4B^{-2}(A^{-1}B)=1$ or $(A^{-1}B)^5 = B^2$. Thus we have
\[
\pi_1(\Sigma) = \langle A, B \, |\, (A^{-1}B)^5 = A^3 = B^2 \rangle\,.
\]

Now imagine that the Borromean rings could be deformed (isotoped is the technical word) to the trivial link. Where would the surgery curves go? Notice that the surgery curve is on the boundary of a solid torus centered on one of the rings. It would still be on the boundary of such a solid torus after the deformation. The surgery curve is parallel to the link component in this solid torus. It would still be parallel after deformation. It also links the component once. After the deformation, it would still have to link once.

In fact one way to define how many times a loop links a simple ring is to look at the class of the loop in the fundamental group of the complement. The fundamental group of the complement of a single ring with no crossings has just one generator and no relations. Call the generator $F$. The only words that may be made in this group are $F^n$ so it is isomorphic to the integers, and we say a loop homotopic to $F^n$ links the ring $n$ times.

The tail of the loop represented by the surgery curve in Figure \ref{bor} might get deformed into something $CRAZY$, so the resulting loop that would have to be killed in the fundamental group of the complement of the trivial $3$-component link would have the form $(CRAZY)P^{\pm 1}(CRAZY)^{-1}$ and this will be trivial exactly when $P^{\pm 1}=1$, i.e., $P=1$.

Thus if the Borromean rings could be deformed into the trivial $3$-component link, the fundamental group of the  Poincar\'e homology sphere would have to be
\[
\langle P, Q, R \, |\, P=Q=R=1\rangle\,.
\]
In other words, it would be the trivial group.

Let's see, we can abelianize the  fundamental group of the Poincar\'e homology sphere. This means we are assuming that $XY=YX$, etc.. Using this in the last three relations in $\pi_1(\Sigma)$ as expressed via $X$, $Y$, and $Z$ generators, gives $X^{-1}=1$, $Y^{-1}=1$ and $Z^{-1}=1$ and this implies that the abelianization of the fundamental group of the Poincar\'e homology sphere is trivial. Maybe we should assume that this means that this space is just a $3$-sphere???  We wouldn't be the first to make this guess. Perhaps it is possible to unlink the Borromean rings after all. Hmm, our yappy friend is reminding us of a leash that represented a trivial loop in an abelianization, but was not trivial in the fundamental group. We need to think.

Well, the fundamental group of the  Poincar\'e sphere does not look trivial. Is this good enough? Consider a different group presentation. The following group is trivial, but it takes work to prove it. (Try.)
\[
\langle C, D \, |\, C^{-1}D^2C=D^3,\ D^{-1}C^2D=C^3\rangle\,.
\]

To see that the fundamental group of the deleted Poincar\'e homology sphere is not trivial, it would suffice to construct a surjective  homomorphism from it to some other non-trivial group. We could do so now and end the paper. However, sometimes it is entertaining when the old man on the porch starts telling stories. This reminds us of something related to an icosahedron, so we'll tell a few more stories.


\section*{The icosahedron}
\begin{figure}[!ht]  
\centering
\begin{subfigure}[b]{0.45\textwidth}
\includegraphics[width=\textwidth]{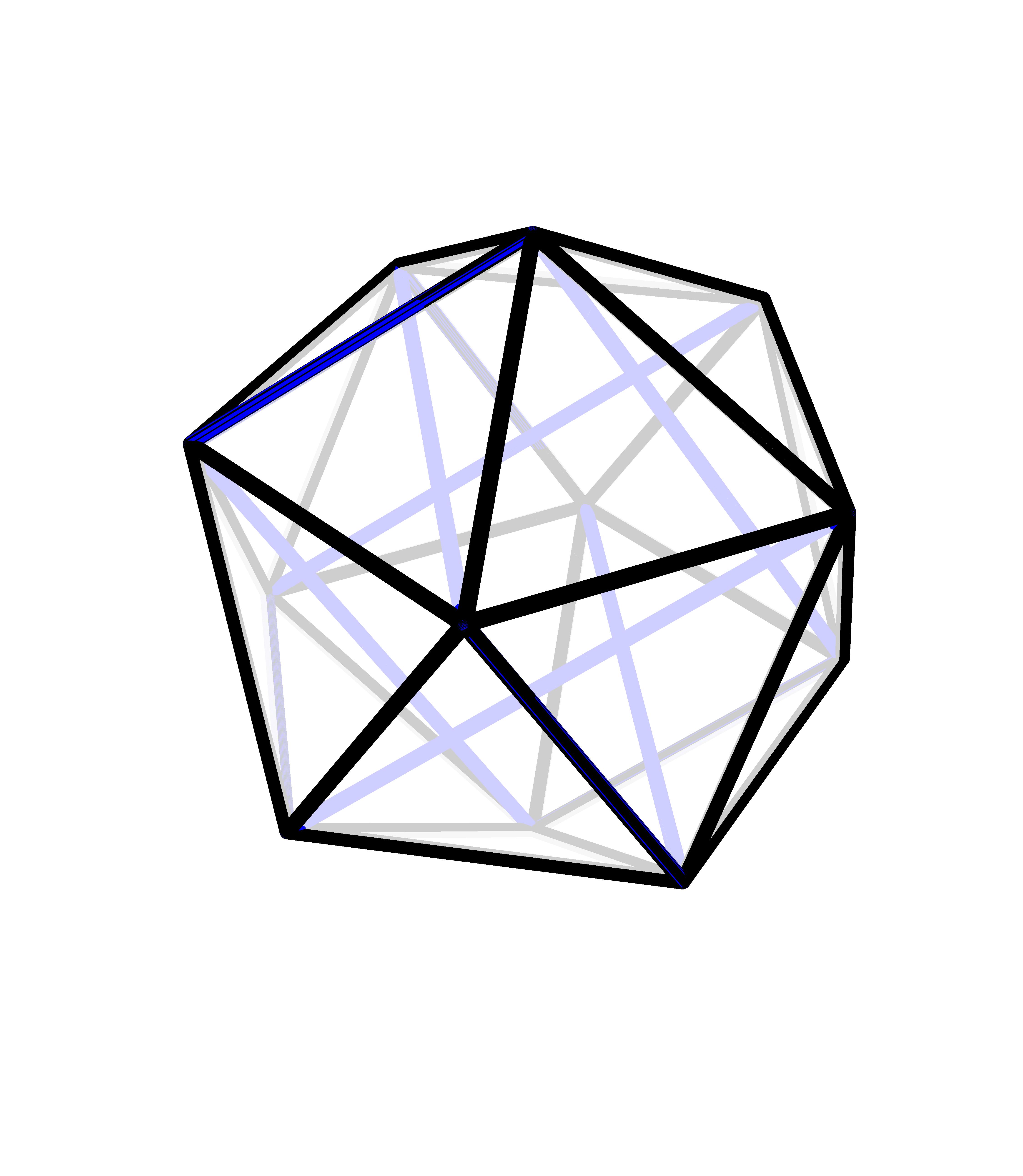}
\caption{Rectangles in an Icosahedron}\label{10a}
\end{subfigure}
\begin{subfigure}[b]{0.45\textwidth}
\includegraphics[width=\textwidth]{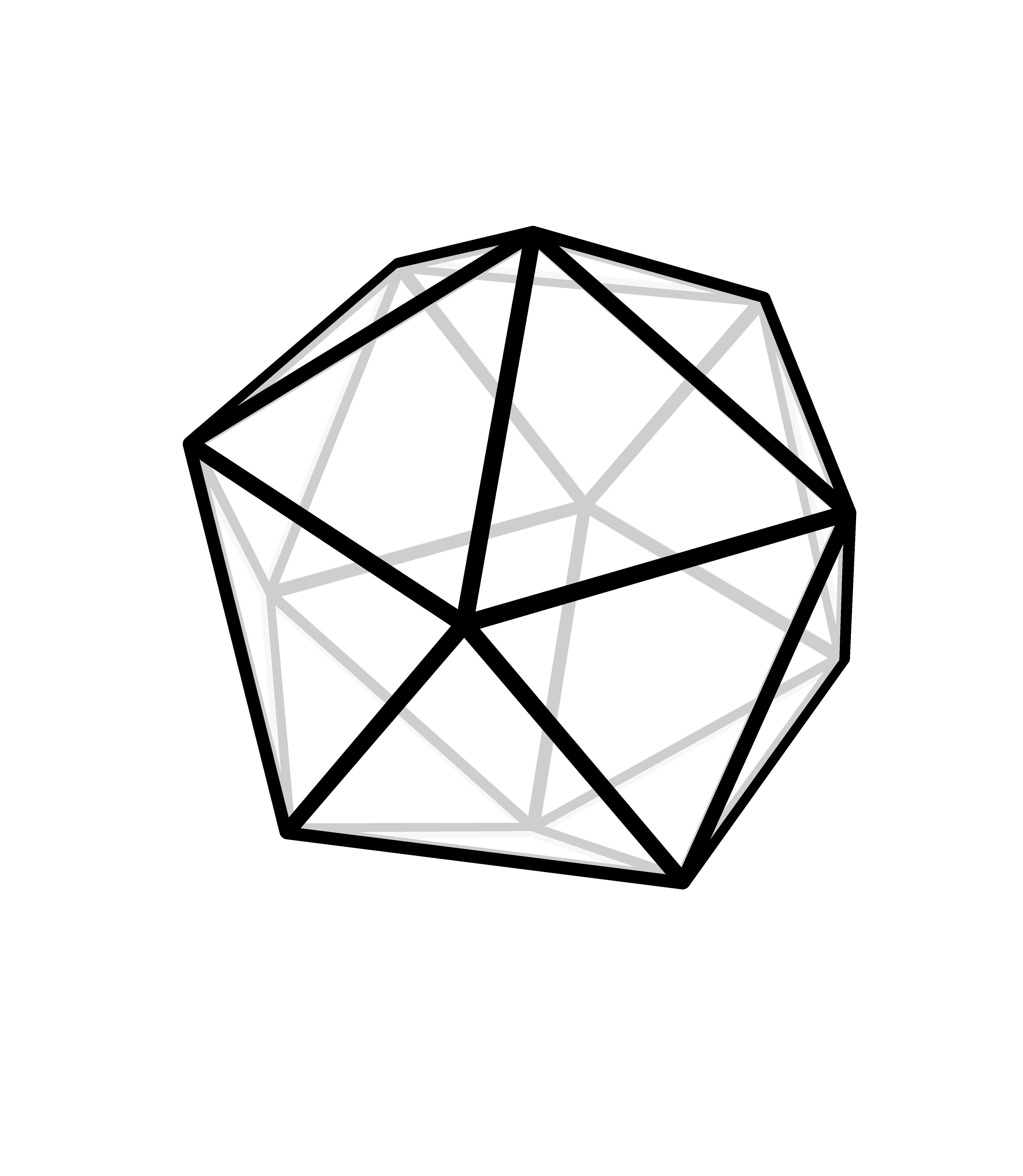}
\caption{Icosahedron}\label{10b}
\end{subfigure}
\caption{Icosahedron around three rectangles}\label{icosa}
\end{figure}

Now return to the three linked rectangles form Figure \ref{3rect} and add $34$-inch struts connecting the corner of each rectangle to the four corners of the other rectangles that are closest to it. The result is the structure displayed in Figure \ref{10a}. If one removes the long edges of each original rectangle, the remaining figure is called an icosahedron. It is the structure represented by the MAA logo, and it is displayed in Figure \ref{10b}.

The icosahedron is one of the five Platonic solids. These solids have captured the imaginations of people for many years. One early model of the orbits of the planets was based on placing one Platonic solid inside of the next. It started by placing the octahedron in the icosahedron. Why don't we try to do the same thing in one way?

\begin{figure}[!ht]  
\centering
\begin{subfigure}[b]{0.45\textwidth}
\includegraphics[width=\textwidth]{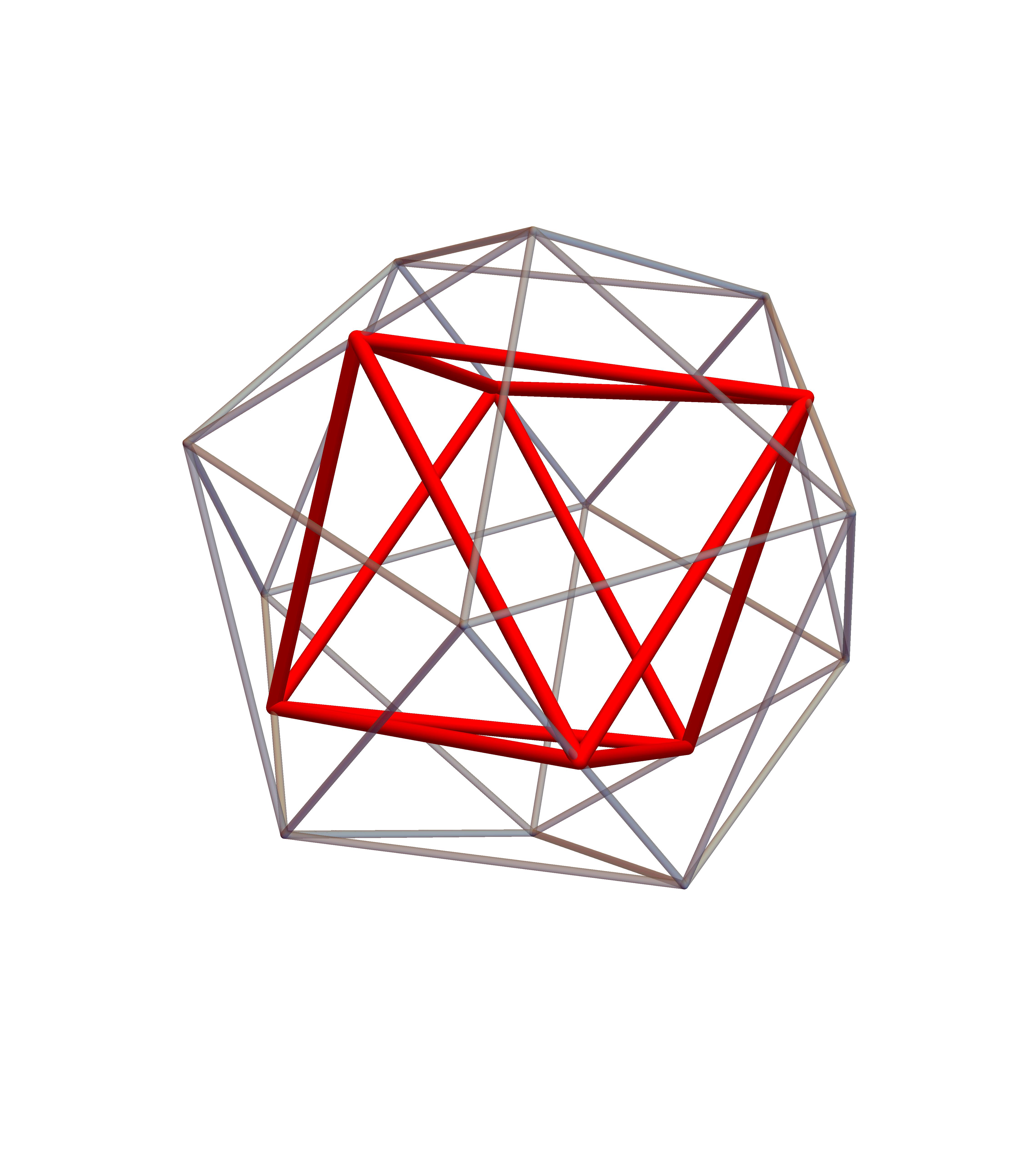}
\caption{Both}\label{11a}
\end{subfigure}
\begin{subfigure}[b]{0.45\textwidth}
\includegraphics[width=\textwidth]{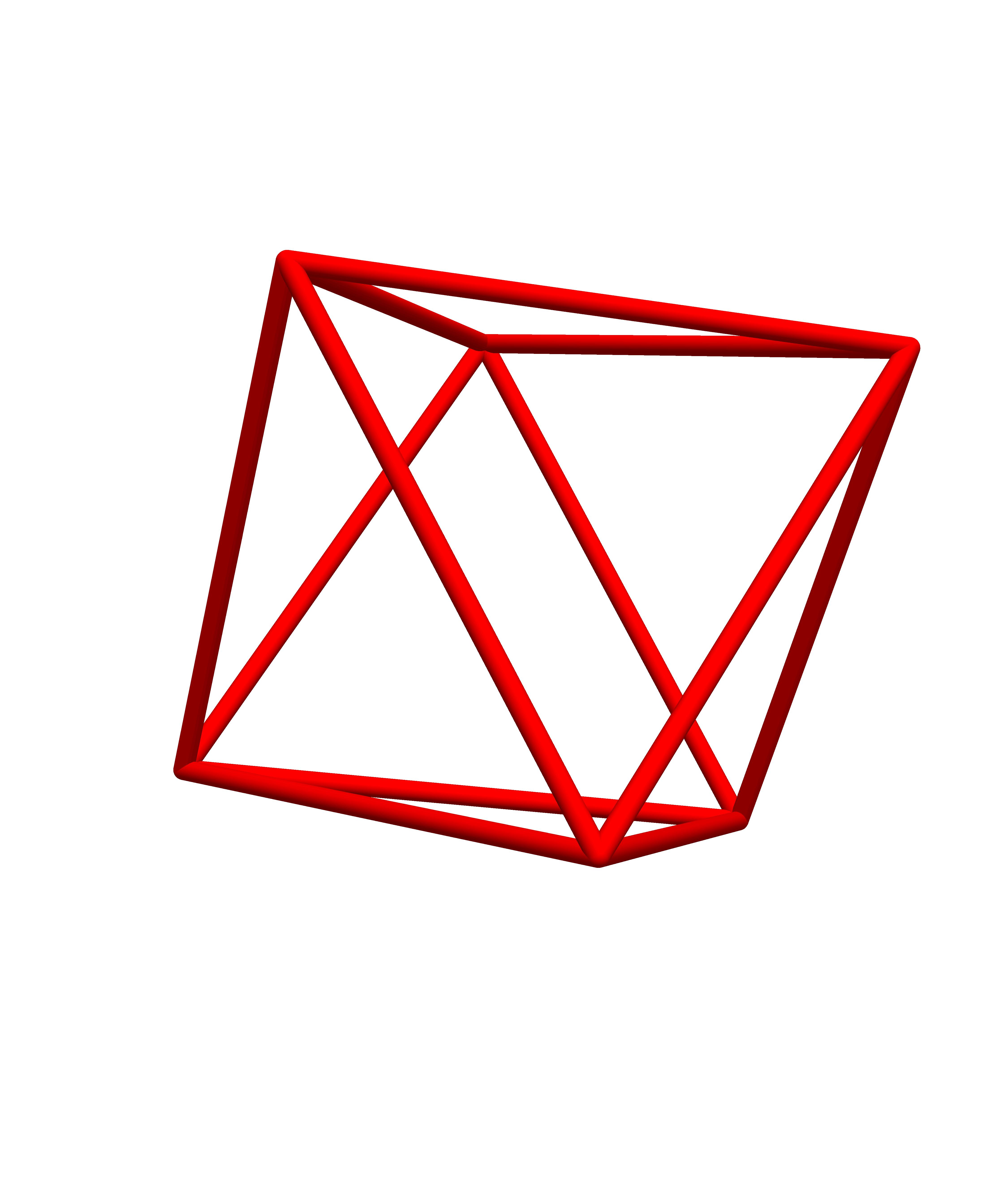}
\caption{Octahedron}\label{11b}
\end{subfigure}
\caption{Octahedron in an Icosahedron}\label{oct-icosa}
\end{figure}

The convex hull of three congruent, mutually perpendicular line segments meeting at their centers is a regular octahedron. The axes of the three original rectangles exactly meet this condition, so we can add the edges of a regular regular octahedron to our figure by connecting the mid points of the short edges of the original three rectangles. This is displayed in Figure \ref{11a} with just the octahedron in Figure \ref{11b}.

This is fun, why stop. We started with three rectangles, and thus had six short edges. We put the six vertices of a red octahedron at the midpoints of these first six edges. We then added more short edges until each rectangle corner met a total of five edge ends. This gave a total of $4$ corners per rectangle times $3$ original rectangles times $5$ (original short edge plus four new short edges) short edge ends. Thus there are $60$ short edge ends and $30$ short edges. It looks like we can fit $5=30/6$ octahedra inside the icosahedron in this way. Let's try with
\begin{enumerate}
\item a red octahedron,
\item an orange octahedron,
\item a yellow octahedron,
\item a green octahedron, and
\item a blue octahedron.
\end{enumerate}


It works. The result is an icosahedral compound of octahedra as displayed in Figure~\ref{comp}.

\begin{figure}[!ht]  
\hskip30bp
\includegraphics[width=70mm]{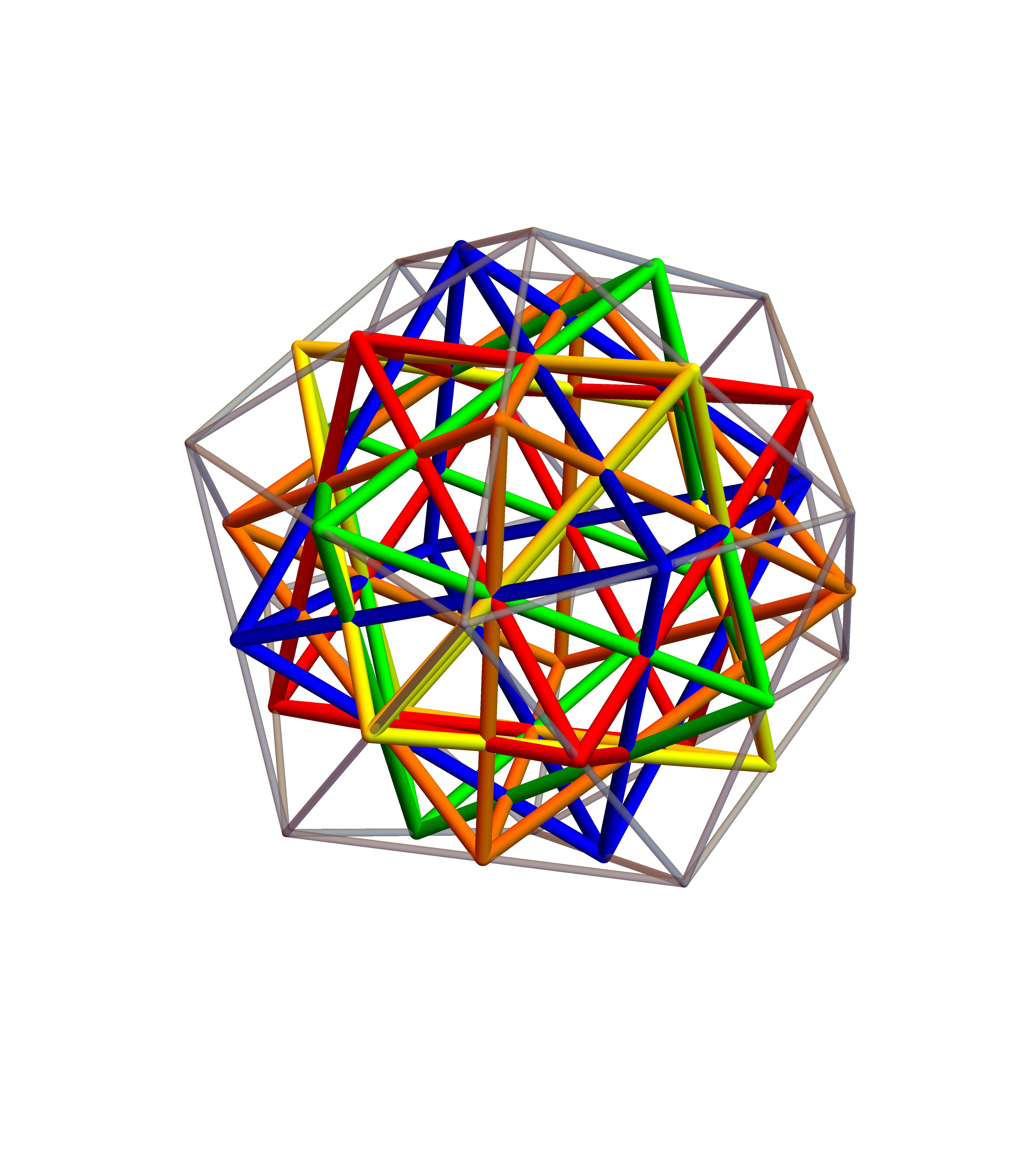}\hskip30bp\includegraphics[width=70mm]{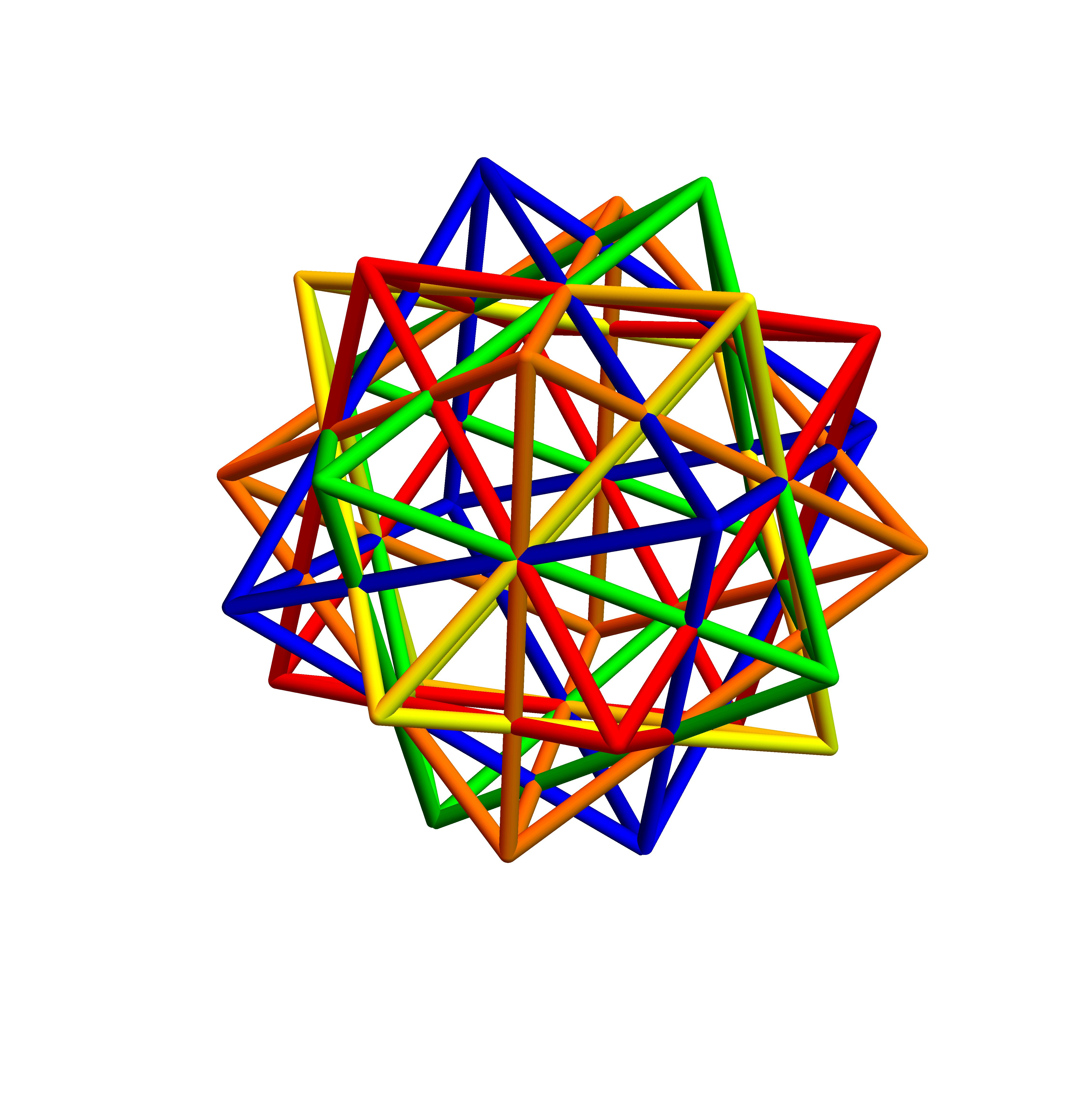}
\caption{Compound of Octahedra in and out of an Icosahedron}\label{comp}
\end{figure}

Consider the orientation-preserving symmetries of the icosahedron. These are all rotations, since they must fix the center of the icosahedron. There is a $1/5\,\text{th}$ right hand rotation about the directed line connecting
the bottom vertex of the icosahedron to the top vertex. This rotation takes: 1 -- the red octahedron to 2 -- the orange octahedron, which moves to the location of 3 -- the yellow octahedron, which moves to 4 -- the green octahedron, which moves to 5 -- the blue octahedron, which moves to 1 -- the red octahedron. We will summarize this by labeling this rotation $(12345)$.

There is also a $1/3\,\text{rd}$ right-hand rotation about the directed line connecting the centroid of the bottom front triangle to the centroid of the opposite triangle.
This rotation takes 1 -- the red octahedron, to 2 -- the orange octahedron, to 3 -- the yellow octahedron, which moves to 1 -- the red. This should be labeled by $(123)$.

Finally, there is a $1/2$ rotation about the line connecting the top of the
red octahedron to the bottom of the red octahedron. This interchanges the 4 -- green, and the 5 -- blue octahedra, as well as interchanges the 2 -- orange, and the 3 -- yellow octahedra. We label it by $(23)(45)$.

Every rotational symmetry of the icosahedron must be one of these three types because the axis of rotation will meet the surface of the icosahedron in two opposite fixed-points. The only possibilities for one of these surface fixed-points are a vertex, the center of a triangle, or the midpoint of an edge. There are $12$ vertices, so $6$ opposite vertex pairs. Each vertex pair has $1/5$, $2/5$, $3/5$, and $4/5$ rotations for a total of $6\times 4 = 24$ rotations. As we saw in the construction, there are $30$ edges, so $15$ pairs of opposite edges, and $15$ more rotations. Each vertex meets $5$ triangle corners for a total of $5\times 12$ triangle corners and $60/3 = 20$ triangular faces. Each of the $20/2 = 10$ opposite triangle pairs contributes $1/3$ and $2/3$ rotations for a total of $20$ more rotations. Combined with the trivial motion $(1)$, we see that there are a total of $24+15+20+1 = 60$ orientation-preserving symmetries.

It is not difficult to see the group operations in terms of the labels. For example, the label $(153)$ corresponds to a map taking $5$ to $3$ and $3$ to $1$, so it is clear that the inverse and the cube are the maps corresponding to
\[
(153)^{-1} = (135)\, \quad\text{and} \quad (153)^3 = (1)\,.
\]
Similarly,
\[
(135)\circ (12)(34) = (12345)\, \quad \left((12)(34)\right)^2 = (1)\,\quad \text{and} \quad (12345)^5 = (1)\,.
\]
Thus,
\[
\left((153)^{-1}(12)(34)\right)^5 = (153)^3 = \left((12)(34)\right)^2\,.
\]

The set of all functions permuting the numbers $1, ..., 5$ forms a group known as the permutation group on $5$ symbols. It is denoted by $\mathfrak{G}_5$ and has $5! = 120$ elements. The labeling scheme demonstrates that the orientation-preserving symmetry group of the icosahedron is isomorphic to the index-$2$ subgroup of $\mathfrak{G}_5$ known as the alternating group of $5$ symbols, $A_5$.

Oh, this reminds us where we were going. We were going to explain why the fundamental group of the  Poincar\'e homology sphere is non-trivial. We just saw that $a=(153)$ and $b=(12)(34)$ satisfy the same relations as $A$ and $B$ in the simplified presentation of the fundamental group of the  Poincar\'e homology sphere. Recall these were $(A^{-1}B)^5 = A^3 = B^2$. Yet the alternating group on $5$-symbols is not trivial.
In fact this shows that the fundamental group of the Poincar\'e homology sphere surjects onto the symmetry group of the icosahedron.

The Poincar\'e  homology sphere is a remarkably beautiful space that appears in many different contexts. Read about several descriptions of it in {\it Eight Faces of the Poincar\'e Homology Sphere}, \cite{8}. To learn more about computing fundamental groups of knot and link complements, cutting out tubes and gluing back in donuts, and many other foundations of low-dimensional topology, Rolfsen's book is a good place to start, \cite{R}.  

The Author would like to thank the referees and Bob Burckel  for very helpful comments on an earlier draft.

\bibliography{1-2020logo}

\bigskip\noindent
{\bf Biography} Dave Auckly is amazed by how fast his children are growing up and treasures time with his family. He has been actively involved in mathematical outreach for his entire career. He served as the Assoicate Director of MSRI for three years, and is one of the co-founders of the Navajo Natoin Math Circles Project. This paper is based on an activity used at the Navajo Nation Math Circles Summer Camp.

\begin{figure}[!ht]  
\center{
\includegraphics[width=70mm]{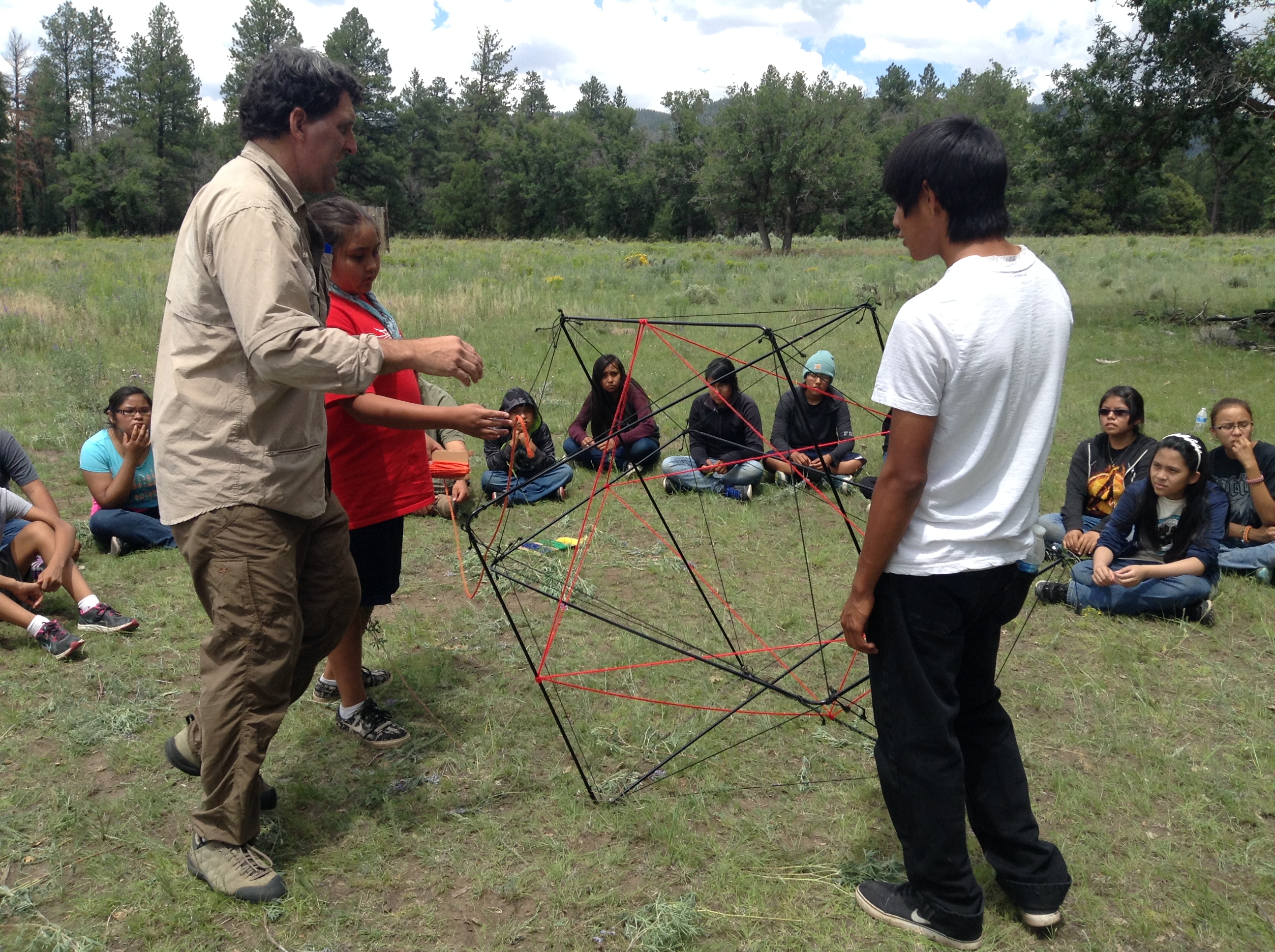}}
\caption{The icosahedral tent at math camp}\label{comp}
\end{figure}

\end{document}